\documentclass[10pt]{amsart}
\usepackage{amsmath,amssymb,amstext,amsthm,amscd,epic,eepic}

\theoremstyle{plain}
\newtheorem{proposition}{Proposition}

\newtheorem{theorem}{Theorem}
\newtheorem*{theorem*}{Theorem}

\newtheorem*{claim*}{Claim}
\newtheorem{conjecture}{Conjecture}
\newtheorem*{conjecture*}{Conjecture}
\newtheorem*{question*}{Question}

\theoremstyle{definition}
\newtheorem{definition}{Definition}
\newtheorem*{definition*}{Definition}

\theoremstyle{remark}
\newtheorem{remark}{Remark}
\newtheorem*{remark*}{Remark}
\newtheorem*{remarks*}{Remarks}

\newtheorem*{acknowledgments}{Acknowledgments}

\newcommand{\QQ}{\mathbb{Q}}

\newcommand{\Mbar}{\overline{\mathcal{M}}}

\def\({\left(}
\def\){\right)}
\def\<{\langle}
\def\>{\rangle}

\def\mr{\mathfrak{r}}
\def\lfun{\longrightarrow}
\def\cM{{\mathcal M}}
\def\ocM{{\overline{\cM}}}

\begin{document}

\title{On independence of generators of the tautological rings}

\author{D.~Arcara}

\author{Y.-P.~Lee}

\address{Department of Mathematics, University of Utah,
Salt Lake City, UT 84112-0090, USA}

\email{arcara@math.utah.edu}

\email{yplee@math.utah.edu}

\begin{abstract}
We prove that all monomials of $\kappa$-classes and $\psi$-classes are 
independent in $R^k(\ocM_{g,n})/R^k(\partial\ocM_{g,n})$ 
for all $k \leq [g/3]$.
We also give a simple argument for $\kappa_l \neq 0$ in $R^l(\mathcal{M}_g)$
for $l \leq g-2$.
\end{abstract}

\maketitle



\section{Introduction}


\subsection{Tautological rings}
Let $\Mbar_{g,n}$ be the moduli stacks of stable curves.
$\Mbar_{g,n}$ are proper, irreducible, smooth Deligne--Mumford stacks.
The Chow rings $A^*(\Mbar_{g,n})$ over $\QQ$ are isomorphic to
the Chow rings of their coarse moduli spaces.
The tautological rings $R^*(\Mbar_{g,n})$ are subrings of
$A^*(\Mbar_{g,n})$, or subrings of $H^{2*}(\Mbar_{g,n})$ via cycle maps,
generated by some ``geometric classes'' which will be described below.

The first type of geometric classes are the \emph{boundary strata}.
$\Mbar_{g,n}$ have natural stratification by topological types.
The second type of geometric classes are the Chern classes of
tautological vector bundles. These includes cotangent classes $\psi_i$,
Hodge classes $\lambda_k$ and $\kappa$-classes $\kappa_l$.

To give a precise definition of the tautological rings,
some natural morphisms between moduli stacks of curves will be used.
The \emph{forgetful morphisms}
\begin{equation} \label{e:ft}
 \operatorname{ft}_i: \Mbar_{g,n+1} \to \Mbar_{g,n}
\end{equation}
forget one of the $n+1$ marked points.
The \emph{gluing morphisms}
\begin{equation} \label{e:gl}
 \Mbar_{g_1,n_1 +1} \times \Mbar_{g_2, n_2 +1} \to \Mbar_{g_1 + g_2, n_1+n_2},
 \quad \Mbar_{g-1,n+2} \to \Mbar_{g,n},
\end{equation}
glue two marked points to form a curve with a new node.
Note that the boundary strata are the images (of the repeated applications)
of the gluing morphisms, up to factors in $\QQ$ due to automorphisms.

\begin{definition}
The system of tautological rings $\{R^*(\Mbar_{g,n})\}_{g,n}$ is the
smallest system of $\QQ$-unital subalgebra (containing classes of
type one and two, and is) closed under the forgetful and gluing morphisms.
\end{definition}

The study of the tautological rings is one of the central problems
in moduli of curves. 
The readers are referred to \cite{rV} and references therein for 
many examples and motivation.

\subsection{Main result}

Let $\cM_{g,n}$ be the moduli stack of smooth $n$-pointed curves.
Let $\partial \Mbar_{g,n} := \Mbar_{g,n} \setminus \cM_{g,n}$.
The main result of this paper is the following theorem.

\begin{theorem} \label{t:1}
The monomials generated by $\kappa$'s and $\psi$'s have no relations in
$R^k(\ocM_{g,n})/R^k(\partial\ocM_{g,n})$ for $k\leq[g/3]$ and all $n$.
\end{theorem}
We note that $R^*(\partial\ocM_{g,n})$ and the quotient 
$R^k(\ocM)/R^k(\partial\ocM)$
are defined in Equation~\eqref{e:3} and the paragraph following it.

\subsection{Motivation: Faber's conjecture}
The formulation of Theorem~\ref{t:1} is motivated by Faber's conjecture on the
structure of the tautological rings.

One of the guiding problems in the study of tautological rings of
moduli of curves is the set of conjectures proposed by C.~Faber, and
R.~Pandharipande, E.~Looijenga \ldots.
(For a survey, the readers are referred to \cite{rV}.) 
In \cite{cF}, Faber conjectures that (Conjecture~1 b.) 
$\kappa_1,\kappa_2,\ldots,\kappa_{[g/3]}$ generate the ring $R^*(\cM_{g})$, 
with no relations in degree $\leq [g/3]$.
\emph{This statement will be referred to as Faber's conjecture in this paper.}
The generation statement has been proved by S.~Morita \cite{sM}, and 
by E.~Ionel \cite{eI}.
Therefore, the remaining part of Faber's conjecture would be the 
independence statement.

An expert in tautological ring can immediately notice the relation between
Faber's conjecture and Mumford's conjecture on the 
\emph{stable} cohomology ring of $\cM_{g,n}$.
In particular, the following theorem of Madsen--Weiss (Mumford's conjecture)
is similar to Faber's conjecture.

\begin{theorem}[\cite{MW} Madsen--Weiss]
\[ 
 H^*(\cM_g) = \mathbb{Q}[\kappa_1, \kappa_2,\ldots],
\]
for $* \le 2g/3$ (the stable range).
\end{theorem}

There is a generalization of Mumford's conjecture to $\cM_{g,n}$, proved
by E.~Looijenga conditional to Mumford's conjecture.

\begin{theorem}[\cite{eL2} Proposition~2.1]
\[
H^*(\cM_{g,n}) = \mathbb{Q}[\kappa_1, \kappa_2,\ldots, \psi_1,\ldots,\psi_n]
\]
in the stable range.
\end{theorem}

Therefore, it is natural to generalize Faber's conjecture accordingly.

\begin{conjecture}\label{c:1}
$\kappa_1, \ldots,\kappa_{[g/3]}, \psi_1, \ldots, \psi_n$ generate the ring 
$R^*(\cM_{g,n})$, with no relations in degree $\leq [g/3]$.
\end{conjecture}
The generation statement follows immediately from Morita and Ionel's results
cited above. What is in question is the independence statement.

Note that there is a natural sequence of tautological rings
\begin{equation} \label{e:3}
R^k(\partial\ocM_{g,n}) \lfun R^k(\ocM_{g,n}) \lfun R^k(\cM_{g,n}) \lfun 0,
\end{equation}
where $R^k(\cM_{g,n})$ is defined to be the restriction of 
$R^k(\Mbar_{g,n})$ and $R^k(\partial\ocM_{g,n})$
is defined to be the pushforward of the normalized boundary divisors 
via gluing morphisms \eqref{e:gl}.
The exactness in the middle of this sequence was conjectured by Faber and 
Pandharipande in \cite{cFrP}, and proved in some cases.
This statement will be referred to as the 
\emph{boundary tautological class conjecture}.
Theorem~\ref{t:1} implies Conjecture~\ref{c:1} if the boundary tautological
class conjecture holds.

\subsection{Invariance constraints}

The main tool employed is the Invariance Constraints,
Theorem~5 of \cite{ypL2}, originally Invariance Conjecture~1 in \cite{ypL1}.
More recently, Faber--Shadrin--Zvonkine, and independently R.~Pandharipande
(and the second author), have given a very simple geometric proof of this
statement. See Section~3 of \cite{FSZ}.

\begin{theorem} \label{t:ic1}
There exist a series of linear operators
\begin{equation} \label{e:inv}
 \mathfrak{r}_l : R^k (\ocM_{g,n}) \to R^{k-l+1} (\ocM^\bullet_{g-1, n+2}),
\end{equation}
where ${}^\bullet$ denotes moduli of possibly disconnected curves.
\end{theorem}
In particular, these operators give an inductive process of finding 
the existence or non-existence of tautological equations. 

The definition of $\mathfrak{r}_l$ is defined by graph operations.
The strata can be conveniently presented by their (dual) graphs, which
can be described as follows.
To each stable curve $C$ with marked points, one can associate a dual graph
$\Gamma$.
Vertices of $\Gamma$ correspond to irreducible components.
They are labeled by their geometric genus.
Assign an edge joining two vertices each time the two components intersect.
To each marked point, one draws an half-edge incident to the vertex,
with the same label as the point.
Now, the stratum corresponding to $\Gamma$ is the closure of the subset of all
stable curves in $\Mbar_{g,n}$ which have the same topological type as $C$.
For each dual graph $\Gamma$, one can decorate the graph by assigning
a monomial, or more generally a polynomial,
of $\psi$ to each half-edge and $\kappa$ classes to each vertex.
The tautological classes in $R^k(\Mbar_{g,n})$ can be
represented by $\QQ$-linear combinations of \emph{decorated graphs}.

Define three graph operations the spaces of decorated graphs
$\{ \Gamma \}$.

\begin{itemize}
\item \textbf{Cutting edges.} Cut one edge and create two new half-edges.
Label two new half-edges with $i,j \notin \{1,2,\ldots,n\} $ 
in two different ways.
Produce a formal sum of 4 graphs by decorating extra $\psi^l$ to $i$-labeled 
new half-edges with coefficient $1/2$ and by decorating extra $\psi^l$ to 
$j$-labeled new half-edges with coefficient $(-1)^l/2$  . 
(By ``extra'' decoration we mean that $\psi^l$ is multiplied by
whatever decorations which are already there.)
Produce more graphs by proceeding to the next edge.
Retain only the stable graphs. 
Take formal sum of these final graphs.

\item \textbf{Genus reduction.}
For each vertex, produce $l$ graphs. Reduce the genus of this given
vertex by one, add two new half-edges.
Label two new half-edges with $i,j$ and decorate them
by $\psi^{l-1-m}, \psi^m$ (respectively) where $0 \le m \le l-1$.
Do this to all vertices, and retain only the stable graphs.
Take formal sum of these graphs with coefficient $\frac{1}{2} (-1)^{m+1}$.

\item \textbf{Splitting vertices.} Split one vertex into two.
Add one new half-edge to each of the two new vertices.
Label them with $i,j$ and decorate them by $\psi^{l-1-m}, \psi^m$
(respectively) where $0 \le m \le l-1$.
Produce new graphs by splitting the genus $g$ between the two new vertices
($g_1, g_2$ such that $g_1+g_2=g$),
and distributing to the two new vertices the (old) half-edges which
belongs to the original chosen vertex, in all possible ways.
The $\kappa$-classes on the given vertex are split between the two 
new vertices in a way similar to the half-edges.
That is, consider each monomial of the $\kappa$-classes 
$\kappa_{l_1} \ldots \kappa_{l_p}$ on the split vertex as labeled by 
$p$ special half-edges. 
When the vertex splits, distribute the $p$ special edges in all possible way.
Do this to all vertices, and retain only the stable graphs.
Take formal sum of these graphs with coefficient $\frac{1}{2} (-1)^{m+1}$.
\end{itemize}

\begin{definition}
$\mathfrak{r}_l (\Gamma)$ is defined to be the formal sum of the outputs of
the above three operations.
\end{definition}
These operations actually descend to $R^k(\Mbar_{g,n})$.
That is, if $E=\sum c_m \Gamma_m = 0 \in R^k(\Mbar_{g,n})$ is a tautological
equation, then
\[
 \mathfrak{r}_l (E) =0.
\]
This is the content of Theorem~\ref{t:ic1}.

\begin{remarks*}
(i) In fact, only the $l=1$ case will be used in this paper.

(ii)
The image of $\mathfrak{r}_l$ lies in the connected components of 
$\ocM^\bullet_{g-1, n+2}$ whose curves have at most two disconnected 
components.

(iii) The definition of these operators are inspired by Givental's study of 
deformation of (axiomatic) Gromov--Witten theories \cite{aG}.
The interested readers are referred to \cite{ypL1} and \cite{ypL2} for
details.
\end{remarks*}

\begin{acknowledgments}
We wish to thank R.~Cavalieri, C.~Faber, R.~Pandharipande, and R.~Vakil for 
helpful discussions. 

The second author is partially supported by NSF and 
an AMS Centennial Fellowship.
\end{acknowledgments}

\section{The Proof}

Theorem~\ref{t:1} will be proved by induction on $(g,n)$, in the 
lexicographic order.

The case $g \leq 2$ is obvious.
Assume now the statement holds for all genera up to $g-1$ and for all $n$.
The case $(g,0)$ will be first proved via Theorem~\ref{t:ic1}.
The following proposition will be used in the proof.

\begin{proposition} \label{p:1}
$\kappa_l \neq 0$ in $R^*(\mathcal{M}_g)$ for all $l \leq g-2$.
\end{proposition}

\begin{proof}
Recall that a nodal genus $g$ curve is said to have rational tails if
one of the irreducible component is smooth of genus $g$,
and $\mathcal{M}_{g,n}^{rt}$ is the moduli stack of genus $g$ nodal curves
with rational tails.
It was shown in \cite{GJV} Theorem~2.5 that
\begin{equation} \label{e:GJV}
 \pi_* (\psi_1^{l+1} \psi_2^{g-l-1}) = c \kappa_{g-2},
\end{equation}
where $\pi$ is the forgetful morphism 
$\mathcal{M}_{g,2}^{rt} \to \mathcal{M}_g$ and
\[
 c = \frac{(2g-1)!!}{(2l+1)!! (2g-2l-3)!!}.
\]
On the other hand, it is well-known that
\[
 \pi_* (\psi_1^{l+1} \psi_2^{g-l-1}) = \kappa_l \kappa_{g-l-2} + \kappa_{g-2}.
\]
See, for example, Lemma~2 of \cite{ypL2}.
Therefore,
\[
 (c-1) \kappa_{g-2} = \kappa_l \kappa_{g-l-2}.
\]
since $c-1 \neq 0$, and $\kappa_{g-2} \neq 0$ (Theorem~2 in \cite{cF}),
the product $\kappa_l \kappa_{g-l-2}$ is non-zero.
This implies that each of the two factors is non-zero.
\end{proof}

\begin{remark}
(i) We suspect that Proposition~\ref{p:1} is known to some experts, but
we were not able to locate a reference. 
In fact, experts we have consulted with were not sure about its status.

(ii) As shown above, it is an immediate consequence of Faber--Looijenga's
\emph{Socle Theorem} \cite{eL1} \cite{cF} plus the fact $c \neq 1$.
It might be possible to prove $c \neq 1$ without the full power of \cite{GJV}.

(iii) Equation \eqref{e:GJV} is part of Faber's \emph{Intersection Number
Conjecture}, first established by E.~Getzler and R.~Pandharipande in \cite{GP}
conditional to Virasoro conjecture of $\mathbb{P}^2$, which was later
established by A.~Givental \cite{aG}.
I.~Goulden, D.~Jackson, and R.~Vakil recently give an alternative proof for 
$n \le 3$ \cite{GJV}.
\end{remark}

\subsection{Case $n=0$}

Assume that
\[
  E = \sum_{I} c_{I} \kappa^I + \sum_m c_m \Gamma_m =0
\]
is a tautological equation in $R^k(\ocM_g)$ for $k \leq [g/3]$.
Here $\kappa^I$ are (distinct) monomials of $\kappa$-classes, 
and $\Gamma_m$ are in the image of $R^k(\partial \Mbar_{g,n})$ via
\eqref{e:3}.
The goal is to show that $c_I=0$ for all $I$.

By Theorem~\ref{t:ic1}, we have that $\mr_1(E)=0$ in
$R^k(\ocM^\bullet_{g-1,2})$.
Let $\kappa^I$ be a monomial in $\kappa$'s of degree $k$.
Now we will analyze the output $\mathfrak{r}_1(\kappa^I)$.
It is easy to see that the following term appears in 
$\mathfrak{r}_1(\kappa^I)$ (splitting the vertex) and does not appear in
$\mathfrak{r}_1(\kappa^J)$ for $J \ne I$ or $\mathfrak{r}_1(\Gamma_m)$:
\begin{eqnarray*}
\setlength{\unitlength}{0.01cm}
\begin{picture}(200,100)(0,0)
\thicklines
\put(40,50){\circle*{20}}
\put(0,70){$g$-$1$}
\put(10,0){$\kappa^I$}
\put(40,50){\line(1,0){70}}
\put(90,65){$i$}
\put(135,65){$j$}
\put(130,50){\line(1,0){70}}
\put(200,50){\circle*{20}}
\put(190,70){$1$}
\end{picture}
\end{eqnarray*}
where $g-1$, $1$ are the genera and $i,j$ are the new half-edges.

Suppose that $k<[g/3]$.
The combination of the following three facts implies $c_I=0$:
\begin{itemize}
\item It only appears in $\mr_1(\kappa^I)$.
\item All monomials in $\kappa$'s of degree $k$ are independent
in $R^k(\ocM_{g-1,1})$ by induction hypothesis.
\item $\mr_1(E)=0$.
\end{itemize}

Let us now assume that $k=[g/3]$.
If $\kappa^I$ is a monomial in $\kappa$'s of degree $[g/3]$ such that
$\kappa^I \neq \kappa_{[g/3]}$,
one can write $\kappa^I=\kappa^{I_1} \kappa^{I_2}$ with $\kappa^{I_a}$ 
a monomial in $\kappa$'s of degree $d_a>0$ ($a=1,2$) such that 
$d_1+d_2=[g/3]$.
Since all monomials in $\kappa$'s of degree $d_a$ are independent in
$R^{d_a}(\ocM_{3d_a,n})$ ($a=1,2$) for all $n$ by induction, the term
\begin{eqnarray*}
\setlength{\unitlength}{0.01cm}
\begin{picture}(200,100)(0,0)
\thicklines
\put(40,50){\circle*{20}}
\put(0,70){$3d_1$}
\put(20,0){$\kappa^{I_1}$}
\put(40,50){\line(1,0){70}}
\put(90,65){$i$}
\put(135,65){$j$}
\put(130,50){\line(1,0){70}}
\put(200,50){\circle*{20}}
\put(190,70){$3d_2$}
\put(190,0){$\kappa^{I_2}$}
\end{picture}
\end{eqnarray*}
in $\mr_1(\kappa^I)$ is independent from any other term appearing in 
$\mr_1(E)$.
Again, since this term only appears in $\mr_1(\kappa^I)$, and not in 
$\mr_1$ of any
other element of $R^{[g/3]}(\ocM_g)$, the fact that its coefficient in
$\mr_1(E)$ must be zero implies that the coefficient of $\kappa^I$ in $E$ 
is also $0$.

The last case to consider is the coefficient of $\kappa_{[g/3]}$.
Suppose that it is nonzero, i.e., suppose that
\[
 E=a\kappa_{[g/3]}+\sum_m c_m\Gamma_m=0
\]
is a tautological equation in $R^{[g/3]}(\ocM_g)$ with $a \neq 0$.
Then, taking the image of this equation in $R^{[g/3]}(\cM_g)$, we would 
obtain that $\kappa_{[g/3]}=0$, which is not true by Proposition~\ref{p:1}
as $g\geq3$ .
Therefore the coefficient of $\kappa_{[g/3]}$ must also be zero.
This concludes the proof of the $(g,0)$ case.

\subsection{Case $n=1$}

Let
\[
  E = \sum_{I, J} c_{IJ} \kappa^I \psi_1^J + \sum_m c_m \Gamma_m =0
\]
be a tautological equation in $R^k(\Mbar_{g,1})$.
Theorem~\ref{t:ic1} implies that $\mr_1(E)=0$ in $R^k(\ocM^\bullet_{g-1,2})$.
The goal is to show that all $c_{IJ} =0$.

Following the same technique, one easily conclude $c_{IJ}=0$ in the following
two cases:
\begin{enumerate}
\item[(a)] $k < [g/3]$.
\item[(b)] $J=0$ and $\kappa^I \neq \kappa_{[g/3]}$.
\end{enumerate}
In both cases, the proof is exactly the same as in $n=0$ case.
Notes that it does not matter where we attach the extra half-edge (marking), 
since the induction hypothesis says the
statement holds for genus $\leq g-1$ and all $n$.

Let us now consider terms of the form 
$\kappa^I \psi_1^J \in R^{[g/3]}(\Mbar_{g,1})$, with $\kappa^I$ a 
monomial of degree $d_I = [g/3]-d_J$.
Suppose that $0<d_J<[g/3]$.
The term
\begin{eqnarray*}
\setlength{\unitlength}{0.01cm}
\begin{picture}(270,100)(0,0)
\thicklines
\put(40,50){\circle*{20}}
\put(20,70){$3d_I$}
\put(20,0){$\kappa^I$}
\put(40,50){\line(1,0){70}}
\put(90,65){$i$}
\put(135,65){$j$}
\put(130,50){\line(1,0){70}}
\put(200,50){\circle*{20}}
\put(180,70){$3d_J$}
\put(200,50){\line(1,0){70}}
\put(280,50){$\psi_1^J$}
\end{picture}
\end{eqnarray*}
in $\mr_1(\kappa^I\psi_1^J)$ is independent from any other 
term appearing in $\mr_1(E)$.
As before, this term only appears in $\mr_1(\kappa^I\psi_1^J)$ (note that it
is important here that $\psi_1^{d_J} \neq \psi_2^{d_J}$ in 
$R^{d_J}(\ocM_{3d_J,2})$, which
follows by induction because $0<d_J<[g/3]$), and therefore the coefficient of
$\kappa^I\psi_1^J$ in $E$ is $0$.

Therefore, we proved that, if $E=0$ is a tautological equation in
$R^k(\ocM_{g,1})$, then all coefficients of monomials in $\kappa$'s and
$\psi_1$ are $0$, except possibly for the coefficients of $\kappa_{[g/3]}$ and
$\psi_1^{[g/3]}$ in $R^{[g/3]}(\ocM_{g,1})$.

Suppose that 
\[
 E=a\kappa_{[g/3]}+b\psi_1^{[g/3]}+\sum_m c_m \Gamma_m=0
\]
is a tautological equation in $R^{[g/3]}(\ocM_{g,1})$, since all other 
coefficients must vanish.
We will perform two operations on $E$ to obtain two linear equations on the 
coefficients of $E$.

\begin{enumerate}
\item[(i)] Multiply $E$ by $\psi_1$ and push-forward of $\psi_1 E$ from
$R^{[g/3]+1}(\ocM_{g,1})$ to $R^{[g/3]}(\ocM_g)$.
The push-forward is easy to perform via the following substitution
\begin{equation} \label{e:kappa}
 \kappa_{[g/3]}=\pi^*\kappa_{[g/3]}+\psi_1^{[g/3]}.
\end{equation}
As $E=0$, the push-forward of $\psi_1 E$ must vanish. Thus
\[
 (2g - 2) a \kappa_{[g/3]} + a \kappa_{[g/3]} + b \kappa_{[g/3]} = 0 
\]
in $R^{[g/3]}(\ocM_g)/R^{[g/3]}(\partial\ocM_g)$.
Since $\kappa_{[g/3]}\neq0$, we obtain that
$$ (2g - 2) a + a + b = 0. $$

\item[(ii)] Push-forward of $E$ from $R^{[g/3]}(\ocM_{g,1})$ to 
$R^{[g/3]-1}(\ocM_g)$.
Equation \eqref{e:kappa} implies
$$ a \kappa_{[g/3]-1} + b \kappa_{[g/3]-1} = 0 $$
in $R^{[g/3]-1}(\ocM_g)/R^{[g/3]-1}(\partial\ocM_g)$.
Since $\kappa_{[g/3]-1}\neq0$, we have that
$$ a + b = 0. $$
\end{enumerate}
(i) and (ii) together imply that $a=b=0$.
This completes the $n=1$ case.

\subsection{Case $n \geq 2$}

Suppose now the statement is proved for $(g, \leq n-1)$ and 
for all genera up to $g-1$.
We will show it holds for $(g,n)$, with $n\geq2$.

Let
\[
  E = \sum_{I, J} c_{IJ} \kappa^I \psi^J + \sum_m c_m \Gamma_m =0
\]
be a tautological equation in $R^k(\Mbar_{g,n})$.
Theorem~\ref{t:ic1} implies that $\mr_1(E)=0$ in
$R^k(\ocM^\bullet_{g-1,n+2})$.
The goal is to show that $c_{IJ}=0$ for all $I,J$.

Again, the case $k<[g/3]$ is very similar to the previous cases and is left
to the readers.

Suppose that $k=[g/3]$.
If $J=0$, the following term only appears in the output 
$\mathfrak{r}_1(\kappa^I)$:
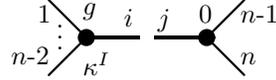
\begin{eqnarray*}
\setlength{\unitlength}{0.01cm}
\begin{picture}(300,100)(-50,0)
\thicklines
\put(-10,100){\line(1,-1){50}}
\put(-10,0){\line(1,1){50}}
\put(-25,70){$1$}
\put(0,35){$\vdots$}
\put(-60,15){$n$-$2$}
\put(40,50){\circle*{20}}
\put(35,80){$g$}
\put(35,0){$\kappa^I$}
\put(40,50){\line(1,0){70}}
\put(90,65){$i$}
\put(135,65){$j$}
\put(130,50){\line(1,0){70}}
\put(200,50){\circle*{20}}
\put(190,70){$0$}
\put(200,50){\line(1,1){50}}
\put(200,50){\line(1,-1){50}}
\put(245,70){$n$-$1$}
\put(245,15){$n$}
\end{picture}
\end{eqnarray*}
It is independent of other class in the output of $\mathfrak{r}_1(E)$
by induction hypothesis.
Therefore, the induction hypothesis implies that $c_{I0}=0$ for all $I$.

If both $I \neq 0$ and $J \neq 0$, the monomial is of the form
$\kappa^I \psi^J$, with $\kappa^I$ degree $d_1$ and $\psi^J$ degree 
$d_2$, such that $d_1+d_2=[g/3]$ and $d_1, d_2 \neq 0$.
The following tautological class only appears in the output of
$\mathfrak{r}_1 (\kappa^I \psi^J)$:
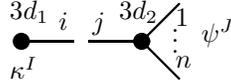
\begin{eqnarray*}
\setlength{\unitlength}{0.01cm}
\begin{picture}(225,100)(25,0)
\thicklines
\put(40,50){\circle*{20}}
\put(25,80){$3d_1$}
\put(25,0){$\kappa^I$}
\put(40,50){\line(1,0){70}}
\put(90,65){$i$}
\put(135,65){$j$}
\put(130,50){\line(1,0){70}}
\put(200,50){\circle*{20}}
\put(170,80){$3d_2$}
\put(280,50){$\psi^J$}
\put(200,50){\line(1,1){50}}
\put(200,50){\line(1,-1){50}}
\put(245,70){$1$}
\put(240,35){$\vdots$}
\put(245,15){$n$}
\end{picture}
\end{eqnarray*}
and is independent of other output classes by induction hypothesis.
Therefore the coefficient $c_{IJ}$ of $\kappa^I \psi^J$ in $E$ is $0$.

The next case is $I=0$. 
If $\psi^J$ is a monomial in $\psi$'s of degree $[g/3]$, and
$\psi^J \neq \psi_l^{[g/3]}$ for any $l$, 
one can write $\psi^J=\psi^{J_1} \psi^{J_2}$, with $\psi^{J_a}$ a monomial in
$\psi$'s of degree $d_a$ ($j=1,2$), $d_1+d_2=[g/3]$, such that
$J_1$ and $J_2$ do not share a common half-edge.
Then the term
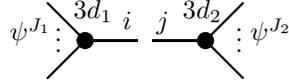
\begin{eqnarray*}
\setlength{\unitlength}{0.01cm}
\begin{picture}(275,100)(-25,0)
\thicklines
\put(-10,100){\line(1,-1){50}}
\put(-10,0){\line(1,1){50}}
\put(0,35){$\vdots$}
\put(40,50){\circle*{20}}
\put(25,80){$3d_1$}
\put(-60,50){$\psi^{J_1}$}
\put(40,50){\line(1,0){70}}
\put(90,65){$i$}
\put(135,65){$j$}
\put(130,50){\line(1,0){70}}
\put(200,50){\circle*{20}}
\put(170,80){$3d_2$}
\put(260,50){$\psi^{J_2}$}
\put(200,50){\line(1,1){50}}
\put(200,50){\line(1,-1){50}}
\put(240,35){$\vdots$}
\end{picture}
\end{eqnarray*}
in $\mr_1(\psi^J)$ implies that the coefficient of $\psi^J$ in $E$ is $0$.

Therefore, we proved that, if $E=0$ is a tautological equation in
$R^k(\ocM_{g,n})$, then all coefficients of monomials in $\kappa$'s and
$\psi$'s are $0$, except possibly for the coefficients of $\psi_1^{[g/3]}$,
$\psi_2^{[g/3]}$, $\ldots$, $\psi_n^{[g/3]}$ in $R^{[g/3]}(\ocM_{g,n})$.

Suppose that 
\[
E=\sum_{l=1}^n a_l \psi_l^{[g/3]} + \sum_m b_m \Gamma_m =0
\]
is a tautological relation in $R^{[g/3]}(\ocM_{g,n})$, where the 
$\Gamma_m$'s are elements of $R^{[g/3]}(\partial\ocM_{g,n})$.

If we multiply by $\psi_n$ and then push-forward to $R^{[g/3]}(\ocM_{g,n-1})$,
we obtain that
$$ \sum_{l=1}^{n-1} (2g - 2 + n - 1) a_l \psi_l^{[g/3]} + a_n \kappa_{[g/3]} =
0 $$
in $R^{[g/3]}(\ocM_{g,n-1})/R^{[g/3]}(\partial\ocM_{g,n-1})$.
Since all monomials in $\kappa$'s and $\psi$'s of degree $[g/3]$ are
independent in $R^{[g/3]}(\ocM_{g,n-1})/R^{[g/3]}(\partial\ocM_{g,n-1})$ by
induction hypothesis, one has
$$ a_1 = a_2 = \cdots = a_n = 0. $$
This completes the proof of Theorem~\ref{t:1}.

\end{document}